\documentstyle{article}
\input amssym.def
\input amssym.tex
\input BoxedEPS
\SetepsfEPSFSpecial

\begin{document}
\title{On the p-affine surface area}
\author{Mathieu Meyer \and Elisabeth Werner
\thanks{ partially supported by a
grant from the National Science Foundation.
MSC classification   52
}}

\date{}
\maketitle
$$\mbox{Mathieu Meyer}$$
$$\mbox{Universit\'{e} de Marne-la-Vall\'{e}e}$$
$$\mbox{Equipe d'Analyse}$$
$$\mbox{93166 Noisy-le-Grand Cedex, France}$$

\vskip 3mm
$$\mbox{Elisabeth Werner}$$
$$\mbox{Department of Mathematics}$$
$$\mbox{Case Western Reserve University}$$
$$\mbox{Cleveland, Ohio 44106}$$
$$\mbox{and}$$
$$\mbox{Universit\'{e} de Lille 1}$$
$$\mbox{Ufr de Mathematique}$$
$$\mbox{59655 Villeneuve d'Ascq, France}$$

\newpage

\vskip 10mm
There are many links between Differential Geometry and Convexity Theory.
An example of such a link is the affine surface area which has attracted
increased
attention in recent years.
\par
Originally a concept of Affine Differential Geometry, it was introduced by
Blaschke [Bl]
for convex bodies in ${\bf R}^3$ with sufficiently smooth boundary and
extended by
Leichtweiss [L 1] to convex bodies in ${\bf R}^n$ with sufficiently smooth
boundary.
Its definition involves the Gauss curvature of the boundary points of a
convex body.
\par
This explains why it has become important in Convexity Theory: it provides
a tool to ``measure" the boundary structure of a convex body.
This is one of the reasons for the renewed interest in the affine surface
area in recent
years. In many applications, so for instance in the approximation of convex
bodies
by polytopes, one needs to have information about the boundary structure of a
convex body. Therefore it is not surprising that the affine surface area
occurs
naturally in many such approximation results -random as well as non-random-
of the last few years by e.g. B${\mbox{\'a}}$r${\mbox{\'a}}$ny [Ba 1], [Ba
2], Gruber [G 1],
[G 2], Gruber and Glasauer [G-G], Ludwig [Lud], Schneider  [Sch] and
Sch\"{u}tt [S 2], to name only a
few.
\par
Another reason is that various isoperimetric inequalities involving the
affine surface area
are very closely related to other important affine isoperimetric
inequalities (e.g., the
curvature image inequality, the Blaschke-Santal${\mbox{\'o}}$ inequality,
and Petty 's
geominimal surface area inequality) (see e.g., [Lu 1], [P 1], [P 2]).
As an application, it has been proved by Lutwak and Oliker  [Lu-O] that
some of these
inequalities lead to a priori estimates  for a certain class of non-linear
PDE's.
\par
From the  point of view of Convexity Theory as well as the applications it
is a drawback to
have the affine surface area only defined for convex bodies with
sufficiently smooth boundary. To find extension of the affine surface area to
arbitrary convex bodies without any smoothness assumptions on the boundary
had been a
problem that was open for a long time.
Fortunately within the last decade several such extensions
to arbitrary convex bodies have been given, namely by  Leichtweiss [L 2],
Lutwak
[Lu 2], Schmuckenschl\"{a}ger [Schm] (for convex symmetric bodies),
Sch\"{u}tt \& Werner  [S-W]
and Werner [W 1]. These extensions all  coincide as was shown by Dolzmann
and Hug  [D-H] (for the
ones given by Leichtweiss and Lutwak),  Schmuckenschl\"{a}ger (for the ones
given by
Schmuckenschl\"{a}ger and Sch\"{u}tt \& Werner), Sch\"{u}tt  [S
1] (for the ones  given by Leichtweiss and Sch\"{u}tt \& Werner) and Werner
[W 1] (for the ones
given by Sch\"{u}tt
\& Werner and Werner).
\par
In [M-W] we investigated a new class of convex bodies which Lutwak called the
Santal${\mbox{\'o}}$-bodies because of their connection with the
Blaschke-Santal${\mbox{\'o}}$ inequality. It came as a surprise to us that
these bodies
provide yet another (completely different) extension of the affine surface
area to arbitrary
convex bodies which also coincides with the existing ones.
All these extensions have a common feature which illustrates nicely the
fact that the affine
surface area is a link between Convexity Theory and Differential Geometry:
geometric
features of the affine surface area are used to give the desired
extensions. In this way we
also obtain a geometric characterisation of the affine surface area.
From these geometric characterisations we get further insight into the
nature of the affine
surface area and thus the boundary structure of a convex body.  In a recent
paper this idea has
been used in [W 2] to give completely general geometric constructions for
the affine surface
area for arbitrary convex bodies which as special cases give all the (so
far) known definitions.
\par
Lutwak [Lu 3] introduced  a generalization of the affine surface area, the
$p$-affine surface
area.   For $p=1$, the p-affine surface area is just the affine surface
area. Lutwak also showed
in [Lu 3]
that the p-affine surface area satisfies extensions of the known
inequalities involving affine
surface area.  As one of the big varieties of results the following
p-extension of the
affine isoperimetric inequality may be stated:
\par
among  all convex bodies in ${\bf R}^n$
with fixed volume the affine surface area is maximal if and only if the
body is an
ellipsoid.
\par
Therefore we expect that the p-affine surface area will turn out to be an
equally useful
tool as the affine surface area.
\par
Hug [H 1] gave new definitions of the p-affine surface area. He also proved
that this new definitons
give the same p-affine surface area as that defined by Lutwak. (In
particular, Hug generalized the
work of Dolzman-Hug from the p=1 case to arbitrary p). Hug showed that for
the case p=n, the
p-affine surface area is the well-known centro-affine surface area. Thus
the notion of p-affine
surface area connects two important affine geometric functionals.
\par
The purpose of this paper is to give a geometric interpretation of the
p-affine surface area
comparable to the ones given for the affine surface area.
This is done in terms of the generalized Santal${\mbox{\'o}}$-bodies which
are of interest in
their own right.
In the first part of the paper we provide the necessary background and
definitions
and introduce the  generalized Santal${\mbox{\'o}}$-bodies and study some
of their properties.
In the second part of the paper we give geometric interpretations of the
p-affine surface area.
\par
The  analytic expression of p-affine surface area in [Lu 3] for convex
bodies with
positive continueous curvature function makes sense not only for positive
p, but also for
$-n < p \leq 0$ and our geometric interpretaion holds for those p too.
We also introduce a definiton of the p-affine surface area for $p=-n$
together with
its geometric interpretation.
\vskip 10mm
 The authors wish to thank MSRI for the hospitality and the
organizers of the special semester in Convex Geometry  and Geometric Functional
Analysis at MSRI for inviting them. It was during the stay there that part
of the
paper was prepared.  We also want to thank E. Lutwak for his many helpful
comments.

\newpage

Let $K$ be a convex body in ${\bf R}^n$ and $t \in {\bf R}$, $t > 0$.
In [M-W] the Santal${\mbox{\'o}}$-bodies $S(K,t)$ were defined as
$$S(K,t) =\{x\in K : \frac{|K||K^{x}|}{v_n^2} \leq t\},$$
where $|K|$ denotes the n-dimensional volume of the convex body $K$ and
$v_n$ is the volume of the n-dimensional Euclidean unit ball $B(0,1)$.
\vskip 3mm
For a convex body $K$ with sufficiently smooth boundary
the affine surface area $O_1(K)$ is
$$O_1(K)= \int_{S^{n-1}} f_K(u)^{\frac{n}{n+1}} d\sigma (u) =
\int_{\partial K} \kappa (x) ^{\frac{1}{n+1}} d\mu_K (x),$$

where $f_K(u)$ is the Gauss curvature function, that is the reciprocal
of the Gauss curvature $\kappa (x)$ at this point $x \in \partial K$ that
has $u$
as outer normal. $\mu_K$ is the usual surface measure on the boundary
$\partial K$
of $K$ and
$\sigma$ is the spherical Lebesgue measure.
\par
The connection between $O_1(K)$ and the Santal${\mbox{\'o}}$-bodies
is as follows
\begin{equation}
\mbox{lim}_{t\rightarrow \infty} t^{\frac{2}{n+1}}(|K| -|S(K,t)|)
=
\frac{1}{2}
(\frac{|K|}{v_n })^{\frac{2}{n+1}}\hspace{.1in}O_1 (K),
\end{equation}
and thus the left hand side provides an extension of the affine
surface area to arbitrary convex bodies without any smoothness assumptions
on the
boundary of $K$.  The one given by (1)
coincides with the ones given earlier by  [L 2], [Lu 2],  [S-W] and [W 1].
\par
In [Lu 3] Lutwak introduced the p-affine surface area $O_p(K)$.
For a convex body $K$
in
${\bf R}^n$ with positive continuous curvature function it can be written as
$$O_p(K)=\int_{S^{n-1}}
\frac{f_K(u)^{\frac{n}{n+p}}}{h_K(u)^{\frac{n(p-1)}{n+p}}}
d\sigma (u) =
\int_{\partial K} \frac{\kappa (x) ^{\frac{p}{n+p}}}
{<x,N(x)>^{\frac{n(p-1)}{n+p}}} d\mu_K (x),$$
where $h_K$ is the support function of $K$ and $N(x)$ is the outer normal
in $x \in
\partial K$.
\vskip 10mm
By $||.||$ we denote the standard
Euclidean  norm on ${\bf R}^n$, $<^.,^.>$ is the usual inner product on
${\bf R}^n$. $B(a,r)$ is the n-dimensional
Euclidean ball with radius $r$
centered at $a$.
\par
For two sets $A$ and $B$ in ${\bf R}^n$, $\mbox{co}[A,B] = \{\lambda a +
(1-\lambda)
 b: a \in A, b \in B, 0 \leq \lambda \leq 1 \}$ is the convex hull of $A$
and $B$.
\par
Unless stated otherwise we will always assume that a convex body $K$
in ${\bf R}^n$ has its Santal${\mbox{\'o}}$-point
at the origin. Then $0$ is the center of mass of the polar body $K^0$
which may be written as
$$\int_{K^0}<x,y> dy = 0
\hspace{.1in}\mbox{for every}\hspace{.1in} x \in {\bf R}^n.$$
Let $\mbox{int} (K)$ be the interior of $K$ and  for $x\in \mbox{int} (K)$,
$K^x = (K-x)^0=
\{y\in
{\bf R}^n: <y,z-x> \leq 1
\hspace{.1in}\mbox{for all z} \in K \}$ is the polar body of $K$ with respect
to $x$;
$K^0$
denotes the polar body with respect to the Santal${\mbox{\'o}}$-point. Moreover
for
$u \in S^{n-1}$ we will denote by $g_{K,u}(s)$
the $(n-1)$ - dimensional volume of the sections of $K$ orthogonal to u,
that is
$$ g_{K,u}(s) =  |\{z \in K: <z,u>=s\}|.$$

\vskip 20mm
Let $\phi : (-1,1) \rightarrow [0, \infty)$ be a continuous function such that
$$\mbox{lim}_{s\rightarrow 1} \phi(s) = \infty.$$
For $x \in \mbox{int} (K)$ we put
\begin{equation}
\Phi_K(x)= \int_{K^0}\phi(<x,y>)dy.
\end{equation}

If $u \in S^{n-1}$ and $\lambda \in {\bf R}$ are such that $0 \leq \lambda <
\frac{1}{h_{K^0}(u)}$, then $x=\lambda u\hspace{0.05in} \in \mbox{int}(K)$ and
(2) can be written as
\begin{equation}
\Phi_K(x)= \int_{-h_{K^0}(-u)}^{h_{K^0}(u)}g_{K^0,u}(s)\phi(\lambda
s)ds.
\end{equation}

\vskip 5mm

\medskip\noindent {\bf Remark 1}
\par
\begin{it}
The motivation to introduce $\Phi_{K}$ comes from [M-W] where it
was observed that for
$u \in S^{n-1}$ and $\lambda \in {\bf R}$ such that $0 \leq \lambda <
\frac{1}{h_{K^0}(u)}$ we have for $x=\lambda u$

$$|K^x|=
\int_{K^0}\frac{dy}{(1-<x,y>)^{n+1}}=
\int_{-h_{K^0}(-u)}^{h_{K^0}(u)}\frac{g_{K^0,u}(s)}{(1-\lambda s)^{n+1}} ds.
$$
The above expressions (2) and (3) generalize this.
\end{it}
\vskip 5mm
We will eventually be interested in more specific functions $\phi$.
But first let us state some general Lemmas.

\newpage

\medskip\noindent {\bf Lemma 2}
\par
\begin{it}
Let K be a convex body in ${\bf R}^n$ and $\phi$ and $\Phi_K$ as above.
\newline
If $\phi$ is convex,
then $\Phi_K$ is a convex function on $\mbox{int}(K)$.
\newline
If $\phi$ is strictly convex, then so is $\Phi_K$.

\end{it}
\vskip 3mm
\medskip\noindent {\bf Proof }

\par
Let $x_1$ and $x_2$ be in $\mbox{int}(K)$, $0 \leq \lambda \leq 1$.
Then
$$\Phi_K (\lambda x_1 +(1-\lambda) \hspace {0.05in} x_2) = \int_{K^0} \phi(<
\lambda x_1 +(1-\lambda) \hspace{0.05in}x_2,y>)\hspace{0.05in}dy,$$
which, by convexity of $\phi$ is
$$\leq
\int_{K^0} \biggl(\lambda\phi(<
x_1,y>) +(1-\lambda) \phi(<x_2,y>)\biggr)\hspace{0.05in}dy $$
$$= \lambda \Phi_K ( x_1) +
(1-\lambda) \Phi_K (
x_2).$$
If $\phi$ is strictly convex, then $\Phi_K$ is strictly convex, as for
$x_1$ and $x_2$  in $\mbox{int}(K)$, $x_1 \neq x_2$, the n-dimensional
volume of the
set$$\{y\in K^0: <x_1,y>=<x_2,y>\} $$ is equal to 0.

\vskip 10mm
We define now for a function $\phi$ with above properties  and for $t \in
{\bf R},
t>0$
$$S_{\phi}(K,t) = \{x \in K: \Phi_K ( x ) \leq t \}.$$

\vskip 3mm
In the sequel we consider only those $t$ for which  the $S_{\phi}(K,t)$
are non-empty.

\vskip 6mm

Then we have

\medskip\noindent {\bf Lemma 3}
\par
\begin{it}
Let $\phi :(-1,1) \rightarrow [0, \infty)$ be a continuous convex
function such that
\newline
$\mbox{lim}_{s\rightarrow 1} \phi(s) = \infty.$
Then
\vskip 3mm
(i) For all $t>0$, $S_{\phi}(K,t)$ is a convex body.
\newline
If $\phi$ is strictly convex, then so is $S_{\phi}(K,t)$.
\vskip 3mm
(ii)
For every affine transformation A with det A $ \neq 0$, for all $t>0$
$$S_{\phi}(A(K),t) = A(S_{\phi}(K,|\mbox{det}A| t)).$$

\end{it}
\newpage
\medskip\noindent {\bf Proof}

\vskip 3mm
(i) follows immediately from Lemma 2.
\vskip 3mm
(ii) Let $A$  be a affine transformation with $\mbox{det}A \neq 0$. We
can write $A=L+a$, where $L$ is a linear transformation with $\mbox{det}L
\neq 0$
and $a$ is a vector in ${\bf R}^n$. Then, as 0 is the
Santal${\mbox{\'o}}$-point of
$K$, $a$ is the Santal${\mbox{\'o}}$-point of $A(K)$ and thus
$$(A(K))^a
=\{z \in {\bf R}^n: <z,Ay-a> \leq 1 \hspace{0.05in}\mbox{for all}
\hspace{0.05in} y \in K \}$$
$$=\{z \in {\bf R}^n: <L^*z,y> \leq 1 \hspace{0.05in}\mbox{for all}
\hspace{0.05in}y \in K \}$$
$$=\{(L^*)^{-1}w: <w,y> \leq 1 \hspace{0.05in}\mbox{for all}
\hspace{0.05in}y \in K \}$$
$$=
(L^*)^{-1}(K^0).$$
Hence
$$S_{\phi}(A(K),t)=\{z \in A(K): \Phi_{A(K)} ( z ) \leq t \}$$
$$=\{Ax: x \in K, \int_{(A(K))^a} \phi( <Ax-a,y> ) \hspace{0.05in}dy\leq t \}$$
$$=\{Ax: x \in K, \int_{(L^*)^{-1}(K^0)} \phi( <Lx,y> ) \hspace{0.05in}dy\leq t
\}$$
$$=\{Ax: x \in K, |\mbox{det}(L^*)^{-1}|\int_{K^0} \phi( <Lx,(L^*)^{-1}y> )
\hspace{0.05in}dy\leq t
\}$$
$$=A(S_{\phi}(K,|\mbox{det}A|t)).$$

\vskip 10mm
Now we consider special functions $\phi$.
In view of Remark 1 a natural class of functions $\phi$ to consider are
$$\phi_{\beta}(s) = \frac{1}{(1-s)^{\beta}},$$
for $\frac{n+1}{2}\leq \beta $.
\newline
For such $\phi_{\beta}$ we denote $S_{\phi_{\beta}}(K,t)$ by
$S_{\beta}(K,t)$, that
is
$$S_{\beta}(K,t)= \{x \in K:
\int_{K^0}
\frac{dy}{(1-<x,y>)^{\beta}}  \leq t\},$$
or with  $x=\lambda u$, $u \in S^{n-1}$,
$0 \leq \lambda < \frac{1}{h_{K^0}(u)}$,
$$S_{\beta}(K,t)= \{x=\lambda u \in K:
\int_{-h_{K^0}(-u)}^{h_{K^0}(u)} \frac{g_{K^0,u}(s)}
{(1-\lambda s)^{\beta}} ds \leq t\}.$$
\vskip 5mm
In particular for $\beta=n+1$ we get the Santal${\mbox{\'o}}$-bodies
$S(K,t)$ of
[M-W].
\vskip 3mm
In the same way for $y \in \partial K^0$,  $y= \lambda u$, $u \in S^{n-1}$,
$0 \leq \lambda < \frac {1}{h_{K}(u)}$
$$S_{\beta}(K^0,t)=
\{y \in K^0:
\int_{K}
\frac{dx}{(1-<x,y>)^{\beta}}  \leq t\}=$$
$$\{y=\lambda u \in K^0:
\int_{-h_{K}(-u)}^{h_{K}(u)} \frac{g_{K,u}(s)}
{(1-\lambda s)^{\beta}} ds \leq t\}.$$
\vskip 5mm

\medskip\noindent {\bf Remark}
\begin{it}
\par
Instead of the functions $\phi_{\beta}$,  a slightly more general class of
functions $\phi$
can
be considered for which the conclusions of Theorem 6 and Proposition 7 will
still hold;
namely for functions $\phi:(-1,1) \rightarrow [0, \infty)$ that are convex,
continuous and
such that
$$\mbox{lim}_{s \rightarrow 1} (1-s)^\beta \phi(s) =c,$$
where $c$ is a constant.
This is easy to see from Lemma  5 and the proofs of Theorem 6 and
Proposition 7.
\end{it}

\vskip 7mm
As for the Santal${\mbox{\'o}}$-bodies $S(K,t)$ one can give estimates on
the ``size" of
$S_{\beta}(K,t)$ in
terms of ellipsoids. Recall that for a convex body K the Binet ellipsoid $E(K)$
is defined by (see for instance [Mi-P])
$$||u||_{E(K)}^2 = \frac{1}{|K|}\int_K <x,u>^2dx , \hspace{.3in}  \mbox{for
all u}
\in{\bf R}^n.$$
\vskip 5mm
\medskip\noindent {\bf Proposition 4}
\par
\begin{it}
Let K be a convex body in ${\bf R}^n$. Then for $\frac{n+1}{2} \leq \beta
\leq n+1$,

$$d_{n}(t, \beta) \hspace{.1in}E(K^0)
\subseteq S_{\beta}(K,t) \subseteq c_{n}(t, \beta)\hspace{.1in}E(K^0),$$

where $$d_{n}(t, \beta) =\frac{1}{\sqrt3
\hspace{.1in}n}\hspace{.1in}(1-\frac{n |K^0|}{t (\beta -1)})^{\frac{1}{2}}$$

and
$$c_{n}(t, \beta) =\frac{2\sqrt2}{((e-2) \beta (\beta +1))^{\frac{1}{2}}}
(\frac{t}{|K^0|})^{\frac{1}{2}}
(1-\frac{|K^0|}{t})^{\frac{1}{2}}.$$
\vskip 3mm
If $K$ is in addition symmetric, then $c_n(t, \beta)$ can be chosen as follows:
$$c_{n}(t, \beta) =\mbox{min}\{\frac{2\sqrt2}{((e-2) \beta (\beta
+1))^{\frac{1}{2}}}
(\frac{t}{|K^0|})^{\frac{1}{2}}
(1-\frac{|K^0|}{t})^{\frac{1}{2}},
\sqrt2
(1-(\frac{|K^0|}{t})^{\frac{1}{\beta -1}})^{\frac{1}{2}} \}.$$

\end{it}
\vskip 5mm
\medskip\noindent {\bf Remark}
\begin{it}
\par
Observe that for $\frac{n+1}{2} \leq \beta \leq n+1$ $ \frac{c_{n}(t,
\beta)}{d_{n}(t, \beta)}$
is equal to a constant that depends only on $\frac{|K^0|}{t}$. This means
that if
$t \leq c |K^0| $, c a constant,  then $S_{\beta}(K,t)$ has bounded (in
terms of c) Banach Mazur
distance to the ellipsoid $E(K^0)$.

\end{it}
\vskip 3mm
\medskip\noindent {\bf {Proof of Proposition 4}}
\par
To get the right hand side inclusions we proceed exactly as in [M-W],
Theorems 6 and 7.
To get the left hand side inclusion, it is enough to consider the symmetric
case by [F].
Then, as in the proof of Theorem 6 of [M-W],
$$t\leq  g_{K^0,u}(0)\int_0^a
(\frac{(1-\frac{s}{a})^{n-1}}{(1-\lambda s)^{\beta }}
+\frac{(1-\frac{s}{a})^{n-1}}{(1+\lambda s)^{\beta }})ds,$$
where $a=
\frac{n|K^0|}{2g_{K^0,u}(0)}$ and for $u \in S^{n-1}$, $\lambda u \in
\partial S_{\beta}(K,t)$.
\newline
In the case $\lambda a \geq 1$, we get immediately as in [M-W]
$$S_{\beta}(K,t) \supseteq\frac{1}{\sqrt3\hspace{.1in}n}
E(K^0).$$
In the case $\lambda a <1$,
we estimate
$$g_{K^0,u}(0)\int_0^a
(\frac{(1-\frac{s}{a})^{n-1}}{(1-\lambda s)^{\beta }}
+\frac{(1-\frac{s}{a})^{n-1}}{(1+\lambda s)^{\beta }})ds\leq $$
$$g_{K^0,u}(0)\int_0^a
(\frac{(1-\frac{s}{a})^{\beta -2}}{(1-\lambda s)^{\beta }}
+\frac{(1-\frac{s}{a})^{\beta -2}}{(1+\lambda s)^{\beta }})ds=$$
$$\frac{g_{K^0,u}(0) \hspace{.1in}a}{\beta -1}\biggl(\frac{1}{(1-\lambda a)}
+\frac{1}{(1+\lambda a)}\biggr),$$
from which it follows, as in [M-W], that
$$S_{\beta}(K,t) \supseteq\frac{1}{\sqrt3\hspace{.1in}n}
\biggl(1-\frac{n|K^0|}{t(\beta -1)}
\biggr)^\frac{1}{2}
E(K^0).$$
\newpage
To relate the convex bodies $S_{\beta}(K,t)$ and $S_{\beta}(K^0,t)$ to the
p-affine surface
area we need the following Lemma.

\vskip 3mm
\medskip\noindent {\bf Lemma 5}
\par
\begin{it}
(i) Let $\gamma > -1$ and $\beta > \gamma +1$. For $\alpha \in (0,1)$ let

$$I(\alpha)= \frac{\alpha^{\gamma +1}(1-\alpha)^{\beta -(\gamma +1)}}
{2^{\gamma } B(\gamma +1,\beta -(\gamma +1))} \int_0^1
\frac{(1-x^2)^{\gamma} dx}{(1-\alpha x)^{\beta}},$$
where for $x,y > 0$, $B(x,y)= \int_0^1 t^ {x-1}(1-t)^{y-1} dt$ is the
Betafunction.
\newline
Then
$$I(\alpha) \leq 1 \hspace{0.1in}
\mbox{and}\hspace{0.1in} \mbox{lim}_{\alpha \rightarrow 1} I(\alpha) =1.$$

\vskip 3mm
(ii) Let $\gamma > 0$ and for $\alpha \in (0,1)$ let
$$J(\alpha)= \frac{1}{2^{\gamma }\mbox{ln}\frac{1}{1-\alpha}}
\int_{0}^1 \frac{(1-x^2)^{\gamma} dx}{(1-\alpha x)^{\gamma +1}}.$$
Then
$$J(\alpha) \leq 1+ \frac{1}{(\gamma
+1)\mbox{ln}\frac{1}{1-\alpha}}\hspace{0.1in}
\mbox{and}\hspace{0.1in}\mbox{lim}_{\alpha \rightarrow 1} J(\alpha) =1.$$

\vskip 5mm
\end{it}
\medskip\noindent {\bf {Proof}}
\par
(i) Put $x=1-w\frac{1-\alpha}{\alpha}$. Then
$$I(\alpha) =\frac{1}
{2^{\gamma}
B(\gamma +1,\beta -(\gamma +1))}
\int_0^{\frac{\alpha}{1-\alpha}}\frac{w^{\gamma}
(2-\frac{1-\alpha}{\alpha}w)^{\gamma}dw}{(1+w)^{\beta}}.$$
The upper estimate for (i) follows immediately from this
last expression.
\newline
And by the Monotone Convergence Theorem this last expression tends to
$$\frac{1}
{B(\gamma +1,\beta -(\gamma +1))}
\int_0^{\infty}\frac{w^{\gamma}
dw}{(1+w)^{\beta }},$$
which is equal to 1.
\vskip 5mm
(ii) Put $ \alpha = 1- e^{-q}$ and $x= 1- e^{-qs}$. Then
$$  \int_{0}^1 \frac{(1-x^2)^{\gamma} dx}{(1-\alpha x)^{\gamma +1}}
= q \int_{0}^{\infty} \frac{(2-e^{-qs})^{\gamma} ds}
{(1+e^{-q(1-s)}-e^{-q})^{\gamma +1}}.$$
Hence
$$J(\alpha) = \frac{1}{2^\gamma}\int_{0}^{\infty}
\frac{(2-e^{-qs})^{\gamma} ds}
{(1+e^{-q(1-s)}-e^{-q})^{\gamma +1}}$$
$$\leq \int_{0}^{1} \frac{ ds}
{(1+e^{-q(1-s)}-e^{-q})^{\gamma +1}} +
\int_{1}^{\infty} \frac{ds}
{e^{q(s-1)(\gamma +1)}}$$
$$\leq 1+ \frac{1}{q(\gamma +1)}.$$

The Monotone Convergence Theorem again implies that
$\mbox{lim}_{\alpha \rightarrow 1} J(\alpha) =1$.

\vskip 5mm
The following theorem gives a geometric interpretation of the p-affine
surface area
for all $p > -n$.

\vskip 5mm
\medskip\noindent {\bf Theorem 6}
\par
\begin{it}
Let K be a convex body in ${\bf R}^n$ such that $\partial K$ is $C^3$
and has strictly positive Gaussian curvature everywhere. Then
for
$\frac{n+1}{2} < \beta $
$$\mbox{lim}_{t\rightarrow \infty}  (\frac{
t}{c_{n,\beta}})^{\frac{1}{\beta-\frac{n+1}{2}}}(|K| -|S_{\beta}(K,t)|)=
\int_{S^{n-1}}\frac{
f_{K^0}(u)^\frac{1}{2 \beta -(n+1)}}{h_{K^0}(u)^{n-\frac{n+1}{2 \beta -(n+1)}}}
d\sigma(u),$$

where
$c_{n,\beta }=2^\frac{n-1}{2}v_{n-1} B(\frac{n+1}{2}, \beta - \frac{n+1}{2})$.

\vskip 3mm

\end{it}

\vskip 5mm
\medskip\noindent {\bf Remarks}
\par
\begin{it}
(i) Thus we have that
$$\mbox{lim}_{t\rightarrow \infty}  (\frac{
t}{c_{n,\beta}})^{\frac{1}{\beta-\frac{n+1}{2}}}(|K| -|S_{\beta}(K,t)|)=
O_{n(2 \beta -n-2)}(K^0).$$
\vskip 3mm
Especially for $\beta = n+1$ we get
$$O_{n^2 }(K^0)=\mbox{lim}_{t\rightarrow \infty}  (\frac{
t}{c_{n,n+1}})^{\frac{2}{n+1}}(|K| -|S_{n+1}(K,t)|).$$
$S_{n+1}(K,t)$ however is the Santal${\mbox{\'o}}$-body $S(K,t)$ introduced in
[M-W]  and it was shown there that
$$\mbox{lim}_{t\rightarrow \infty}  (\frac{
t}{c_{n,n+1}})^{\frac{2}{n+1}}(|K| -|S(K,t)|)= O_1(K).$$
Thus we get again a special case of a general formula of Hug [H 2] who
showed that for $p > 0$
$$O_p(K)=O_{\frac{n^2}{p}}(K^0).$$
\vskip 3mm
Therefore we also get by Hug's formula and by Theorem 6 that for $p > 0$
$$\mbox{lim}_{t\rightarrow \infty}  (\frac{
t}{c_{n,\beta}})^{\frac{1}{\beta-\frac{n+1}{2}}}(|K| -|S_{\beta}(K,t)|)=
O_{\frac{n}{2 \beta -n-2}}(K).$$

\newpage
(ii) As $(K^0)^0=K$, it follows immediately from Theorem 6 under the same
hypothesis on $K$ that

$$\mbox{lim}_{t\rightarrow \infty}  (\frac{
t}{c_{n,\beta}})^{\frac{1}{\beta-\frac{n+1}{2}}}(|K^0| -|S_{\beta}(K^0,t)|)=
\int_{S^{n-1}}\frac{
f_{K}(u)^\frac{1}{2 \beta -(n+1)}}{h_{K}(u)^{n-\frac{n+1}{2 \beta -(n+1)}}}
d\sigma(u)$$
$$=
O_{n(2 \beta -n-2)}(K).$$
\vskip 5mm
(iii) Using a compactness argument , the method of our proof shows that
the statement of Theorem 6 holds if we only suppose that $\partial K$ is $C^2$
and has strictly positive Gaussian curvature everywhere.
\end{it}
\vskip 7mm
\medskip\noindent {\bf Proof of Theorem 6}
\vskip 3mm
Note first that if $K$ is such that $\partial K$ is $C^3$
and has strictly positive Gaussian curvature everywhere,
then the same holds for $K^0$.
\par
For $u \in S^{n-1}$
let $y \in \partial K^0$ be such that $N(y)=u$.
\par
By assumption the indicatrix of Dupin at $y$ exists and is an ellipsoid.
We can assume that it is a sphere
(see for instance  [S-W]). Let
$\sqrt{\rho}=\sqrt{\rho(u)}$ be the radius of this sphere.

We introduce a coordinate system such that $y =0$ and
$u=N(y) = (0, \ldots 0, -1)$.  $H_0$ is the tangent hyperplane to $\partial
K^0$ in
$y =0$ and $\{ H_s: s \geq 0 \}$ is the family of hyperplanes
parallel to $H_0$
that have non-empty
intersection with $K^0$ and are at distance $s$ from $H_0$.  For $s > 0$,
$H_s^+$ is the
halfspace generated by $H_s$ that contains $y =0$.
For $a\in {\bf R}$, let $z_a=(0,\ldots 0,a)$ and $B_a=B(z_a,a)$ be the
Euclidean ball
with center $z_a$ and radius a.
As in [W], for
$\varepsilon > 0$ there exists  $s_{\varepsilon}=s_{\varepsilon}(u) $ so
that for all
$s \leq s_{\varepsilon}$

$$B_{\rho-\varepsilon} \cap
H_s^+
\subseteq
K^0 \cap H_s^+
\subseteq
B_{\rho+\varepsilon}
\cap H_s^+.$$
We choose $s_0 = \mbox{min} \{s_{\varepsilon}, \frac{\rho -\varepsilon}{2} \}$.
\newline
Define $C_1$ to be the cone tangent to $B_{\rho+\varepsilon}$ at
$H_{s_0} \cap B_{\rho+\varepsilon}$ and choose the minimal $s_1$
so that
$$K^0 \cap \{z: s_0 \leq <z, -u> \leq s_1 \}
\subseteq D=C_1 \cap \{z: s_0 \leq <z, -u> \leq s_1 \}.$$
Then  $K^0$ is contained in the union of the truncated cone $D$  and the cap
$B_{\rho + \varepsilon} \cap H_{s_0}^+ =
\{ z \in B_{\rho + \varepsilon}: <z,-u> \leq s_0 \}$\begin{equation}K^0
\subseteq D \cup (B_{\rho + \varepsilon} \cap H_{s_0}^+).
\end{equation}

Let $P$ be the point of intersection of $\partial K^0$ with the positive
$x_n$-axis.
Let $C_2$ be the spherical cone $C_2 = \mbox{co}[P,B_{\rho - \varepsilon}
\cap H_{s_0}]$
and let $h$ be the height of this cone. Then
\begin{equation}
K^0 \supseteq C_2 \cup (B_{\rho - \varepsilon} \cap H_{s_0}^+).
\end{equation}
\vskip 7mm
Now we want to estimate $g_{K^0,u}(s)$. To do so we switch the origin of the
coordinate system such that the u-coordinate of the centroid of $K^0$ is at
0 and
the positive u-direction coincides with the positive s-direction (see
figure below
for the body containing $K^0$ ).

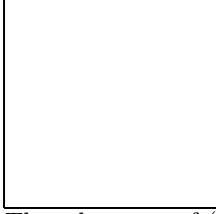
\BoxedEPSF{figure6 scaled 800}
%\begin{figure}
%\epsfbox{figure6}
%\end{figure}

Thus, because of (4)

$$g_{K^0,u}(s) \leq v_{n-1}\biggl[(\rho+ \varepsilon)^2 - \biggl(s-(h_{K^0}(u)
-(\rho+ \varepsilon))\biggr)^2\biggr]^\frac{n-1}{2},$$
if
$$h_{K^0}(u)-s_0 \leq s \leq h_{K^0}(u),$$

and
$$g_{K^0,u}(s) \leq
v_{n-1}\biggl[\frac{s_0(\rho+\varepsilon)+
(\rho+\varepsilon -s_0)(h_{K^0}(u)-s)}{\biggl(2(\rho+\varepsilon)s_0
-s_0^2\biggr)^\frac{1}{2}}\biggr]^{n-1},$$  if
$$-h_{K^0}(-u)\leq s \leq h_{K^0}(u)-s_0.$$
\vskip 3mm

And because of (5)
$$g_{K^0,u}(s) \geq
v_{n-1}\biggl[(\rho- \varepsilon)^2 - \biggl(s-(h_{K^0}(u)
-(\rho- \varepsilon))\biggr)^2\biggr]^\frac{n-1}{2},$$
if
$$h_{K^0}(u)-s_0 \leq s \leq h_{K^0}(u),$$

and
$$g_{K^0,u}(s) \geq
\frac{v_{n-1}}{h^{n-1}}
\biggl[(2s_0(\rho-\varepsilon)-s_0^2)^{\frac{1}{2}}
(s-h_{K^0}(u)+s_0 +h)\biggr]^{n-1},
$$
if
$$h_{K^0}(u)-s_0 -h \leq s \leq h_{K^0}(u)-s_0.$$

\vskip 3mm
Let $\lambda \in {\bf R}, 0 \leq \lambda < \frac{1}{h_{K^0}(u)}$.
Then we have for $x= \lambda u \in \partial S_{\beta}(K,t)$
$$t=\int_{-h_{K^0}(-u)}^{h_{K^0}(u)}
\frac{g_{K^0,u}(s) ds}{(1-\lambda s)^{\beta}}$$

\begin{equation}
\leq v_{n-1}(I_1+I_2),
\end{equation}
where
$$I_1 = \int_{h_{K^0}(u)-s_0}^{h_{K^0}(u)}
\frac{\biggl((\rho+ \varepsilon)^2 - \biggl(s-(h_{K^0}(u)
-(\rho+ \varepsilon))\biggr)^2\biggr)^\frac{n-1}{2} ds}{(1-\lambda
s)^{\beta}}\biggr \}$$
and
$$I_2 = \int_{-h_{K^0}(-u)}^{h_{K^0}(u)-s_0}\frac{
\biggl(s_0(\rho+\varepsilon)+
(\rho+\varepsilon
-s_0)(h_{K^0}(u)-s)\biggr)^{n-1}ds}{\biggl(2(\rho+\varepsilon)s_0
-s_0^2\biggr)^{\frac{n-1}{2}}(1-\lambda s)^{\beta}}.
$$
\vskip 3mm
We consider first $I_1$.
$$I_1=(\rho+ \varepsilon)^{n-1}\int_{h_{K^0}(u)-s_0}^{h_{K^0}(u)}
\frac{\biggl[1 - \biggl(1+\frac{s}{\rho+ \varepsilon}-
\frac{h_{K^0}(u)}{\rho+
\varepsilon} \biggr)^2\biggr]^\frac{n-1}{2} ds}
{(1-\lambda s)^{\beta}}.$$

We put $v=1+\frac{s}{\rho+ \varepsilon}-
\frac{h_{K^0}(u)}{\rho+
\varepsilon}$
and get
$$I_1 \leq \frac{(\rho+ \varepsilon)^n}
{\biggl(1-\lambda (h_{K^0}(u)-(\rho+
\varepsilon))\biggr)^{\beta}}
\int_{0}^{1}
\frac{(1-v^2)^\frac{n-1}{2} dv}
{\biggl(1-\frac{\lambda (\rho+ \varepsilon)v}
{1-\lambda (h_{K^0}(u)-(\rho+
\varepsilon))}\biggr)^{\beta}},$$

which, by Lemma 5 (i) is
$$\leq \frac{2^{\frac{n-1}{2}}(\rho+ \varepsilon)^{\frac{n-1}{2}}
B(\frac{n+1}{2}, \beta - \frac{n+1}{2})}{\lambda ^{\frac{n+1}{2}}
(1-\lambda h_{K^0}(u))^{\beta -\frac{n+1}{2}}}.$$

\vskip 3mm
Using $\rho + \varepsilon \geq s_0$ and $\lambda (s+s_0) \leq \lambda
h_{K^0}(u) < 1$,
we have
$$I_2 \leq$$
$$\frac{\biggl(s_0(\rho+ \varepsilon) +(\rho+
\varepsilon-s_0)(h_{K^0}(u)+h_{K^0}(-u)) \biggr)^n -
\biggl(s_0(\rho+ \varepsilon) +(\rho+
\varepsilon-s_0)s_0 \biggr)^n
}{n(\lambda s_0)^{\beta}(2(\rho+\varepsilon)s_0-s_0^2)^{\frac{n-1}{2}}
(\rho+\varepsilon-s_0)}$$

$$\leq \frac{(\rho+\varepsilon)^{\frac{n+1}{2}}
\biggl(s_0 + h_{K^0}(u)+h_{K^0}(-u)\biggr) ^n}
{n \lambda ^{\beta}(\rho+\varepsilon-s_0) s_0^{\beta+ \frac{n-1}{2}}}.$$
Thus, putting $c_{n,\beta }=2^\frac{n-1}{2}v_{n-1} B(\frac{n+1}{2}, \beta -
\frac{n+1}{2})$, we get

$$t \leq$$
$$\frac{c_{n,\beta }(\rho+
\varepsilon)^{\frac{n-1}{2}}}{\lambda ^{\frac{n+1}{2}}(1-\lambda h_{K^0}(u)
)^{\beta -\frac{n+1}{2}}}
\biggl \{1+\frac{v_{n-1}(\rho+
\varepsilon)(1-\lambda h_{K^0}(u)
)^{\beta -\frac{n-1}{2}}(s_0 + h_{K^0}(u)+ h_{K^0}(-u))^n}
{n c_{n,\beta }(\rho+\varepsilon -s_0)
s_0^{\beta+\frac{n+1}{2}} \lambda^{\beta -\frac{n+1}{2}}}\biggr \}.$$

We choose $\lambda$ so big that
\begin{equation}
\lambda \geq
\frac{1-\mbox{min}\{\varepsilon^{\frac{\beta +1 +\frac{n+1}{2}}{\beta
-\frac{n+1}{2}}},
s_0^{\frac{\beta +1 +\frac{n+1}{2}}{\beta -\frac{n+1}{2}}}
\}}{\mbox{min}\{h_{K^0}(u),\rho-\varepsilon\}}.
\end{equation}
Then
\begin{equation}
t \leq
\frac{c_{n,\beta }(\rho+
\varepsilon)^{\frac{n-1}{2}}h_{K^0}(u) ^{\frac{n+1}{2}}}
{(1-\lambda
h_{K^0}(u) )^{\beta -\frac{n+1}{2}}}(1+ c_1 \varepsilon),
\end{equation}
where $c_1$ is a constant.

\vskip 7mm
On the other hand
$$t \geq v_{n-1}
\int_{h_{K^0}(u)-s_0}^{h_{K^0}(u)}
\frac{\biggl((\rho- \varepsilon)^2 - \biggl(s-(h_{K^0}(u)
-(\rho- \varepsilon))\biggr)^2\biggr)^\frac{n-1}{2} ds}{(1-\lambda
s)^{\beta}},$$
which with $v=1+\frac{s}{\rho- \varepsilon}-
\frac{h_{K^0}(u)}{\rho-
\varepsilon}$ is equal to

$$ \frac{v_{n-1}(\rho- \varepsilon)^n}
{\biggl(1-\lambda (h_{K^0}(u)-(\rho-
\varepsilon))\biggr)^{\beta}} \biggl \{
\int_{0}^{1}
\frac{(1-v^2)^\frac{n-1}{2} dv}
{\biggl(1-\frac{\lambda (\rho- \varepsilon)v}
{1-\lambda (h_{K^0}(u)-(\rho-
\varepsilon))}\biggr)^{\beta}}$$
\begin{equation}
-\int_{0}^{1-\frac{s_0}{\rho -\varepsilon}}
\frac{(1-v^2)^\frac{n-1}{2} dv}
{\biggl(1-\frac{\lambda (\rho- \varepsilon)v}
{1-\lambda (h_{K^0}(u)-(\rho-
\varepsilon))}\biggr)^{\beta}} \biggr \}.
\end{equation}
As (7) holds we get with Lemma 5 (i)
$$t \geq $$
$$\frac{c_{n,\beta }(\rho-
\varepsilon)^{\frac{n-1}{2}}h_{K^0}(u) ^{\frac{n+1}{2}}}
{(1-\lambda
h_{K^0}(u) )^{\beta
-\frac{n+1}{2}}}\biggl(1-\varepsilon-\frac{2^{\frac{n+3}{2}}v_{n-1}
(\rho -\varepsilon)^{\frac{n+1}{2}}(1-\lambda
h_{K^0}(u) )^{\beta -\frac{n+1}{2}}}
{(n+1) c_{n,\beta }
s_0^{\beta} \lambda^{\beta -\frac{n+1}{2}}}\biggr)$$
or, again using (7)
\begin{equation}
t \geq
\frac{c_{n,\beta }(\rho-
\varepsilon)^{\frac{n-1}{2}}h_{K^0}(u) ^{\frac{n+1}{2}}}
{(1-\lambda
h_{K^0}(u) )^{\beta -\frac{n+1}{2}}}(1-c_2 \varepsilon),
\end{equation}
where $c_2$ is a constant.
\vskip 5mm
Thus we get for $x= \lambda u \in \partial  S_{\beta}(K,t)$
$$
\lambda = \frac{1}{h_{S_{\beta}(K,t)^0}(u)}$$
and from (8)
$$\lambda \geq \frac{1}{h_{K^0}(u)}\biggl\{1-\biggl (\frac{c_{n,\beta }(\rho+
\varepsilon)^{\frac{n-1}{2}}h_{K^0}(u) ^{\frac{n+1}{2}}(1+c_1 \varepsilon)}{t}
\biggr )^{\frac{1}{\beta -\frac{n+1}{2}}}\biggr\}$$
respectively from (10)
$$\lambda \leq
\frac{1}{h_{K^0}(u)}
\biggl\{1-\biggl (
 \frac{
c_{n,\beta }(\rho-\varepsilon)^{\frac{n-1}{2}}h_{K^0}(u)
^{\frac{n+1}{2}}(1-c_2 \varepsilon)}{t}
\biggr )^{\frac{1}{\beta -\frac{n+1}{2}}}\biggr\}.$$

\vskip 3mm
Therefore for $t$ big enough
$$(1-d_1 \varepsilon)^{\frac{1}{\beta-\frac{n+1}{2}}}
(\rho (u)-\varepsilon)^{\frac{n-1}{2 \beta-(n+1)}}
h_{K^0}(u)^{-n+\frac{n+1}{2 \beta-(n+1)}} \leq$$

$$\frac{1}{n}(\frac{
t}{c_{n,\beta}})^{\frac{1}{\beta-\frac{n+1}{2}}}\biggl ((\frac{1}
{h_{K^0}(u)})^n - (\frac{1}
{h_{S_{\beta}(K,t)^0}(u)})^n\biggr  ) \leq$$
$$
(1+d_2 \varepsilon)^{\frac{1}{\beta-\frac{n+1}{2}}}
(\rho (u)+\varepsilon)^{\frac{n-1}{2 \beta-(n+1)}}
h_{K^0}(u)^{-n+\frac{n+1}{2 \beta-(n+1)}} ,
$$
where $d_1$ and $d_2$ are constants. This means that for every $u \in S^{n-1}$
$$
\mbox{lim}_{t \rightarrow \infty}
\frac{1}{n}(\frac{
t}{c_{n,\beta}})^{\frac{1}{\beta-\frac{n+1}{2}}}\biggl ((\frac{1}
{h_{K^0}(u)})^n - (\frac{1}
{h_{S_{\beta}(K,t)^0}(u)})^n\biggr  ) = \frac{f_{K^0}(u)^{\frac{1}{2
\beta-(n+1)}}}
{h_{K^0}(u)^{n-\frac{n+1}{2 \beta-(n+1)}}} .
$$
Hence
$$\mbox{lim}_{t \rightarrow \infty}
(\frac{ t}{c_{n,\beta}})^{\frac{1}{\beta-\frac{n+1}{2}}}
\biggl (|K| -|S_{\beta}(K,t)| \biggr  ) =$$
$$\mbox{lim}_{t \rightarrow \infty}\int_{S^{n-1}}
\frac{1}{n}(\frac{
t}{c_{n,\beta}})^{\frac{1}{\beta-\frac{n+1}{2}}}\biggl ((\frac{1}
{h_{K^0}(u)})^n - (\frac{1}
{h_{S_{\beta}(K,t)^0}(u)})^n\biggr  )  d\sigma(u)=$$
$$\int_{S^{n-1}}\mbox{lim}_{t \rightarrow \infty}
\frac{1}{n}(\frac{
t}{c_{n,\beta}})^{\frac{1}{\beta-\frac{n+1}{2}}}\biggl ((\frac{1}
{h_{K^0}(u)})^n - (\frac{1}
{h_{S_{\beta}(K,t)^0}(u)})^n\biggr  ) d\sigma(u)$$
$$=\int_{S^{n-1}}\frac{f_{K^0}(u)^{\frac{1}{2 \beta-(n+1)}}}
{h_{K^0}(u)^{n-\frac{n+1}{2 \beta-(n+1)}}} d\sigma(u) .$$
We still have to justify that we can interchange integration and limit.
This follows  from
Lebesgue's Theorem and the
\newline
\underline{Claim}
\newline
for every $u \in S^{n-1}$
$$t^{\frac{1}{\beta-\frac{n+1}{2}}}\biggl ((\frac{1}
{h_{K^0}(u)})^n - (\frac{1}
{h_{S_{\beta}(K,t)^0}(u)})^n\biggr  )
\leq l(u),$$
where $l$ is a function independent of $t$ and integrable on $S^{n-1}$.
\newline
\underline{Proof of the Claim}
\par
For $u \in S^{n-1}$ let again $ y \in \partial K^0$ be such that $N(y)=u$.
We assume again that the indicatrix of Dupin at $y$ is a
Euclidean ball with radius $\sqrt {\rho}=\sqrt {\rho(u)}$.
We choose $t_0$ so big that $S_{\beta}(K,t_0)$ is a convex body. Then
$S_{\beta}(K,t_0)$
has non-empty interior and there is
 $\alpha > 0$ such
that
$$
B(0,\alpha) \subseteq S_{\beta}(K,t_0) \subseteq
B(0,\frac{1}{\alpha}).
$$
and the same holds for all $t \geq t_0$.
Thus for $x=\lambda u \in \partial S_{\beta}(K,t)$
\begin{equation}
\alpha \leq \lambda=\frac{1}{h_{S_{\beta}(K,t)^0}(u)} \leq \frac{1}{\alpha}.
\end{equation}
As before, we estimate
$$t=\int_{-h_{K^0}(-u)}^{h_{K^0}(u)}
\frac{g_{K^0,u}(s) ds}{(1-\lambda s)^{\beta}}\leq v_{n-1}(I_1+I_2),$$
where $I_1$ and $I_2$ are as above.
From above we get for all $t$
$$I_1 \leq \frac{2^{\frac{n-1}{2}}(\rho+ \varepsilon)^{\frac{n-1}{2}}
B(\frac{n+1}{2}, \beta - \frac{n+1}{2})}{\lambda ^{\frac{n+1}{2}}
(1-\lambda h_{K^0}(u))^{\beta -\frac{n+1}{2}}},$$
which, using (11), can be estimated for all $t\geq t_0$ by
$$\leq k_1 \frac{(\rho+ \varepsilon)^{\frac{n-1}{2}}}{(1-\lambda
h_{K^0}(u))^{\beta
-\frac{n+1}{2}}},$$ where $k_1$ is a constant independent of $t$ and $u$.
Also from above we get for all $t$
$$I_2 \leq
\frac{(\rho+\varepsilon)^{\frac{n+1}{2}}
\biggl(s_0 + h_{K^0}(u)+h_{K^0}(-u)\biggr) ^n}
{n \lambda ^{\beta}(\rho+\varepsilon-s_0) s_0^{\beta+ \frac{n-1}{2}}}.$$
As $s_0 = \mbox{min} \{s_{\varepsilon}, \frac{\rho -\varepsilon}{2} \}
\leq \frac{\rho -\varepsilon}{2}$ and as $K^0$ is bounded, hence contained
in some ball,
for $t \geq t_0$ the last expression can be estimated with (11) by
$$I_2 \leq k_2 \frac{(\rho + \varepsilon)^{\frac{n-1}{2}}}{s_0^{\beta
+\frac{n-1}{2}}},$$
where $k_2$ is a constant.
As $\partial K^0$ is $C^3$and has strictly positive Gaussian curvature
everywhere, $\sigma_\varepsilon=
\mbox{min}_{u\in S^{n-1}}s_\varepsilon(u) > 0$.
Let $R_0 = \mbox{min}_{z \in \partial K^0,1 \leq i \leq n-1}R_i(z)$, where
$R_i(z)$
is the i-th principal radius of curvature at $z \in \partial K^0$. By the
assumptions
on $\partial K$,
$R_0 > 0$  (see [L]) and $\rho \geq R_0$ for all $\rho$.
Thus
$$s_0 \geq \mbox{min}\{\frac{\rho -\varepsilon}{2}, \sigma_\varepsilon\} \geq
\mbox{min}\{\frac{R_0 -\varepsilon}{2}, \sigma_\varepsilon\},$$
which is a stictly positive number independent of $u$ and $t$.
Thus for $t \geq t_0$
$$t \leq k \frac{\rho^{\frac{n-1}{2}}}{(1-\lambda h_{K^0}(u))^{\beta
-\frac{n+1}{2}}},$$ where $k$ is a (new) constant independent of $t$ and $u$.
Therefore
$$\lambda \geq \frac{1}{h_{K^0}(u)}\biggl(1-k
\frac {\rho^{\frac{n-1}{2 \beta-n-1}}}{t^\frac{1}{\beta
-\frac{n+1}{2}}}\biggr)$$
and thus for all $t \geq t_0$
$$t^{\frac{1}{\beta-\frac{n+1}{2}}}\biggl ((\frac{1}
{h_{K^0}(u)})^n - (\frac{1}
{h_{S_{\beta}(K,t)^0}(u)})^n\biggr  ) \leq c \frac{(\rho(u))^{\frac{n-1}{2
\beta-n-1}}}{(h_{K^0}(u))^n },$$ where $c$ is a (new) constant independent
of $t$ and $u$.
\newline
And the function $l(u)= \frac{(\rho(u))^{\frac{n-1}{2
\beta-n-1}}}{(h_{K^0}(u))^n }$ is integrable on $S^{n-1}$.
\vskip 5mm
This proves the Claim and thus the Theorem.

\vskip 7mm
The next Proposition deals with the case $\beta =\frac{n+1}{2}.$
First we want to give a definition for $O_{-n}(K)$ and the motivation for
this definition.
This definition is probably well known though we did not find a reference.
\newline
For $ \lambda \in {\bf R}$, $\lambda \geq 0 $ ([Lu 3], [H 1])
$$O_p(\lambda K)= \lambda ^\frac{n(n-p)}{n+p}O_p(K).$$
Therefore $\tilde{O}_p(K) = O_p(K)^\frac{n+p}{n-p}$ is homogeneous of degree n
and affine invariant.
$$\tilde{O}_p(K) = \biggl(
\int_{S^{n-1}}(\frac{f_{K}(u)^n}{h_{K}(u)^{n(p-1)}})^\frac{1}{n+p} d\sigma(u)
\biggr)^\frac{n+p}{n-p}$$
$$= ||\frac{f_{K}^n}{h_{K}^{n(p-1)}}||_{\frac{1}{n+p}},$$
where for a function $g \in L^q(S^{n-1})$, $||g||_q =(\int_{S^{n-1}} g(u)^q
d\sigma
(u))^\frac{1}{q}$. Thus, if $p \rightarrow -n$,
$$\tilde{O}_p(K) \rightarrow \mbox{max}_{u\in S^{n-1}}
f_{K}(u)^\frac{1}{2}h_{K}(u)^\frac{n+1}{2}= ||
f_{K}^\frac{1}{2}h_{K}^\frac{n+1}{2} ||_\infty.$$
The latter is also an affine invariant and it is therefore natural to put
$$\tilde{O}_{-n}(K)=\mbox{max}_{u\in S^{n-1}}
f_{K}(u)^\frac{1}{2}h_{K}(u)^\frac{n+1}{2}.$$
Proposition 7 respectively Remark (ii) following Proposition 7 gives a
geometric
interpretation of this affine invariant.
\vskip 5mm
\medskip\noindent {\bf Proposition 7}
\par
\begin{it}
Let K be a convex body in ${\bf R}^n$ such that $\partial K$ is $C^3$
and has strictly positive Gaussian curvature everywhere. Then

$$n\mbox{lim}_{t \rightarrow \infty}
\frac{ |K| -|S_{\frac{n+1}{2}}(K,t)| }
{\int_{S^{n-1}}\frac{f_{K^0}(u)^{\frac{1}{n-1}}}{h_{K^0}(u)^{n+1}}
\hspace{0.1in}
\mbox{exp}\biggl(\frac{-t}{2^ {\frac{n-1}{2}}h_{K^0}(u)^{\frac{n+1}{2}}
f_{K^0}(u)^{\frac{1}{2}}}\biggr) d\sigma(u)} = 1.$$

\end{it}
\newpage
\medskip\noindent {\bf Remarks}
\par
\begin{it}
(i) As $(K^0)^0=K$, it follows immediately from Proposition 7 that
$$n\mbox{lim}_{t \rightarrow \infty}
\frac{ |K^0| -|S_{\frac{n+1}{2}}(K^0,t)| }
{\int_{S^{n-1}}\frac{f_{K}(u)^{\frac{1}{n-1}}}{h_{K}(u)^{n+1}} \hspace{0.1in}
\mbox{exp}\biggl(\frac{-t}{2^ {\frac{n-1}{2}}h_{K}(u)^{\frac{n+1}{2}}
f_{K}(u)^{\frac{1}{2}}}\biggr) d\sigma(u)} = 1.$$
\vskip 3mm
(ii) It also follows from Proposition 7 that
$$\mbox{lim}_{t \rightarrow \infty}\frac{1}{t}\mbox{ln}
\biggl(\frac{1}{|K| -|S_{\frac{n+1}{2}}(K,t)| }\biggr)=
\frac{2^{-\frac{n-1}{2}}}{\mbox{max}_{u \in S^{n-1}}(h_{K^0}(u)^{\frac{n+1}{2}}
f_{K^0}(u)^{\frac{1}{2}})}.$$
\end{it}
\par
Indeed, if we put
$B(t)=\int_{S^{n-1}}g(u)\biggl(exp(-\frac{1}{\Phi(u)})\biggr)^td\sigma (u)$,
where $g(u)=\frac{f_{K^0}(u)^{\frac{1}{n-1}}}{h_{K^0}(u)^{n+1}}$ and
$\Phi(u)=2^{\frac{n-1}{2}}f_{K^0}(u)^{\frac{1}{2}}h_{K^0}(u)^\frac{n+1}{2}$
and $A(t)=n(|K| -|S_{\frac{n+1}{2}}(K,t)|)$,
then by Proposition 7, $\frac{A(t)}{B(t)} \rightarrow 1$, as
$t \rightarrow \infty.$ Hence
$$\mbox{ln}\biggl( \frac{A(t)}{B(t)}\biggr)=\mbox{ln}(A(t))-\mbox{ln}(B(t))=
\mbox{ln}(A(t))\biggl(1-\frac{\mbox{ln}(B(t))}{\mbox{ln}(A(t))}\biggr)
\rightarrow 0.$$
As $A(t) \rightarrow 0$ and $B(t) \rightarrow 0$
as $t \rightarrow \infty$, it follows that $\frac
{\mbox{ln}(B(t))}{\mbox{ln}(A(t))}
\rightarrow 1$ as $t \rightarrow \infty$ and thus
$$\mbox{ln}(B(t))^\frac{1}{t} \frac{t}{\mbox{ln}(A(t))}=
\frac{\mbox{ln}(B(t))}{t}\frac{t}{\mbox{ln}(A(t))} \rightarrow1.$$
But
$$(B(t))^\frac{1}{t}=\biggl(\int_{S^{n-1}}g(u)
\biggl(exp(-\frac{1}{\Phi(u)})\biggr)^td\sigma (u)\biggr)^\frac{1}{t}$$
$$= ||exp(-\frac{1}{\Phi})||_{L^t(S^{n-1},gd\sigma)}
\rightarrow ||exp(-\frac{1}{\Phi})||_{L^\infty (S^{n-1})},$$
as $t \rightarrow \infty$. From this (ii) follows.

\vskip 5mm
\medskip\noindent {\bf Proof of Proposition 7}
\vskip 3mm
We use the same notations as in the proof of Theorem 6.
In fact up to the estimates (6) and (9) for $t$ both proofs are identical.
Then we apply Lemma 5 (ii) instead of Lemma 5 (i) to estimate $I_1$
and  get again with
 $v=1+\frac{s}{\rho+ \varepsilon}-
\frac{h_{K^0}(u)}{\rho+
\varepsilon}$
$$I_1 \leq \frac{(\rho+ \varepsilon)^n}
{\biggl(1-\lambda (h_{K^0}(u)-(\rho+
\varepsilon))\biggr)^{\frac{n+1}{2}}}
\int_{0}^{1}
\frac{(1-v^2)^\frac{n-1}{2} dv}
{\biggl(1-\frac{\lambda (\rho+ \varepsilon)v}
{1-\lambda (h_{K^0}(u)-(\rho+
\varepsilon))}\biggr)^{\frac{n+1}{2}}}$$
$$\leq \frac{2^\frac{n-1}{2}(\rho+
\varepsilon)^\frac{n-1}{2}}{\lambda^\frac{n+1}{2}}
\mbox{ln}(1+\frac{\lambda (\rho+ \varepsilon)}{1-\lambda h_{K^0}(u)})
\biggl(1+\frac{2}{(n+1)\mbox{ln}(1+\frac{\lambda (\rho+ \varepsilon)}{1-\lambda
h_{K^0}(u)})}
\biggr).$$

\vskip 3mm
We use the same estimate as in the proof of Theorem 6 for $I_2$ and get with
$\beta=\frac{n+1}{2}$
$$I_2 \leq \frac{(\rho+ \varepsilon)^\frac{n+1}{2}(s_0 +h_{K^0}(u)
+h_{K^0}(-u))^n}
{n(\rho+ \varepsilon-s_0)
\lambda^\frac{n+1}{2}s_0^n}.$$
Thus for
$$
\lambda \geq
\frac{\mbox{max}\{1-\frac{\varepsilon}{2},
1-\frac{e^{\frac{-1}{\varepsilon s_0^n}}}{2}\}}
{\mbox{min}\{h_{K^0}(u),\rho-\varepsilon\}}
$$
$$t \leq (1+ c \varepsilon)v_{n-1}2^\frac{n-1}{2}(\rho
+\varepsilon)^\frac{n-1}{2}
h_{K^0}(u)^\frac{n+1}{2}
\mbox{ln}(1+\frac{\lambda (\rho + \varepsilon)}{1-\lambda h_{K^0}(u)}),$$
where c is a constant.
This implies that
$$\lambda = \frac{1}{h_{(S_{\frac{n+1}{2}}(K,t))^0}(u)} $$
$$\geq \frac{1}{h_{K^0}(u)}\biggl[ 1-\frac{\rho + \varepsilon}{h_{K^0}(u)}
\biggl(\mbox{exp}(\frac{t}
{(1+ c \varepsilon)v_{n-1}2^\frac{n-1}{2}(\rho +\varepsilon)^\frac{n-1}{2}
h_{K^0}(u)^\frac{n+1}{2}})
-1\biggr)^{-1}\biggr].$$

Consequently
$$|K|- |S_{\frac{n+1}{2}}(K,t)| = \frac{1}{n}\int_{S^{n-1} }\biggl
( \frac{1}{h_{K^0}(u)^n} - \frac{1}{h_{(S_{\frac{n+1}{2}}(K,t))^0}(u)^n}
\biggr) d \sigma (u)$$
$$ \leq $$
$$\int_{S^{n-1}} \bigg\{\frac{\rho(u) +\varepsilon}{h_{K^0}(u)^{n+1}}
\mbox{exp} \biggl(\frac{-t}{(1+ c \varepsilon)v_{n-1}2^\frac{n-1}{2}(\rho
+\varepsilon)^\frac{n-1}{2} h_{K^0}(u)^\frac{n+1}{2}}\biggr) $$
$$ \biggl( 1- \frac{1}{
\mbox{exp}(\frac{t}{(1+ c \varepsilon)v_{n-1}2^\frac{n-1}{2}(\rho
+\varepsilon)^\frac{n-1}{2}
h_{K^0}(u)^\frac{n+1}{2}})}\biggr)^{-1} \bigg\}d \sigma (u).$$

For t sufficiently big, this gives the estimate from above.
\vskip 3mm
In the same way
$$t \geq
(1- d\varepsilon)v_{n-1}2^\frac{n-1}{2}(\rho -\varepsilon)^\frac{n-1}{2}
h_{K^0}(u)^\frac{n+1}{2}
\mbox{ln}(1+\frac{\lambda (\rho - \varepsilon)}{1-\lambda h_{K^0}(u)}),$$
where d is a constant.
Therefore
$$\lambda = \frac{1}{h_{(S_{\frac{n+1}{2}}(K,t))^0}(u)} \leq$$
$$\frac{1}{h_{K^0}(u)}\biggl[ 1-\frac{(1 - \varepsilon)(\rho -
\varepsilon)}{h_{K^0}(u)}
\mbox{exp}\biggl(\frac{-t}
{(1- d \varepsilon)v_{n-1}2^\frac{n-1}{2}(\rho -\varepsilon)^\frac{n-1}{2}
h_{K^0}(u)^\frac{n+1}{2}}\biggr)
\biggr]$$
and thus for t sufficiently big
$$|K|- |S_{\frac{n+1}{2}}(K,t)|  \geq$$

$$\int_{S^{n-1}} \frac{(1 - \varepsilon)^2(\rho(u)
-\varepsilon)}{h_{K^0}(u)^{n+1}}
\mbox{exp} \biggl(\frac{-t}{(1-d  \varepsilon)v_{n-1}2^\frac{n-1}{2}(\rho
-\varepsilon)^\frac{n-1}{2} h_{K^0}(u)^\frac{n+1}{2}}\biggr)
 d \sigma (u).$$

\vskip 10mm
 The floating body can also be used to give a
geometric interpretation of the p-affine surface area for certain p.
\par
Recall that for $\delta > 0$, $\delta$ small enough, $K_{\delta}$
is said to be a (convex) floating body of $K$, if it is the intersection of all
halfspaces whose defining hyperplanes cut off a set of volume
$\delta$ of K ([S-W]). More precisely, for $ u \in S^{n-1}$ and for
$ 0 < \delta$ let $a_{\delta}^u$ be defined by
$$|\{x \in K: <x,u> \geq a_{\delta}^u\}| = \delta.$$ Then $K_{\delta}
= \cap_{u \in S^{n-1}}\{x\in K: <x,u> \leq a_{\delta}^u\}$.
Observe that one has always $h_{K_{\delta}}(u)  \leq a_{\delta}(u)$,
with generally strict inequality. There is equality for every $u$ and
every $\delta \leq \frac{|K|}{2}$ whenever $K$ is centrally symmetric
(see [M-R]).
\vskip 3mm
Then we have
\vskip 3mm
\medskip\noindent {\bf Theorem 8}
\par
\begin{it}
Let K be a convex body in ${\bf R}^n$ such that $\partial K$ is $C^2$
and has strictly positive Gaussian curvature everywhere. Then
$$\mbox{lim}_{\delta \rightarrow 0}
c_{n}\frac{|(K_{\delta |K|})^0| -|K^0|}{(\delta |K|)^{\frac{2}{n+1}}}=
\int_{S^{n-1}} \frac{d\sigma(u)}{
f_{K}(u)^{\frac{1}{n+1}}h_{K}(u)^{n+1}},$$

where
$c_{n}=\frac{2 v_{n-1}}{ n+1}$.
\end{it}
\vskip 5mm
\medskip\noindent {\bf Remark}
\par
\begin{it}
Thus
$$\mbox{lim}_{\delta \rightarrow 0}
c_{n}\frac{|(K_{\delta |K|})^0| -|K^0|}{(\delta |K|)^{\frac{2}{n+1}}}
=O_{-n(n+2)}(K).$$
\end{it}
\newpage
\medskip\noindent {\bf Proof of Theorem 8}
\par

For $u \in S^{n-1}$ let $ x \in \partial K$ be such that $N(x)=u$.
By assumption the
indicatrix of Dupin at $x$ exists and is an ellipsoid.
We can again assume that
it is  a Euclidean ball with radius $\sqrt \rho= \sqrt {\rho (u)}$.
With the notations and the coordinate system introduced in  the
proof of Theorem 6 (with x instead of y), for $\varepsilon >0$
there exists $s_0$ such that
$$B_{\rho-\varepsilon} \cap
H_{s_0}^+
\subseteq
K^0 \cap H_{s_0}^+
\subseteq
B_{\rho+\varepsilon}
\cap H_{s_0}^+.$$
We choose $\delta $ so small that
$s_0 > h_{K_{\delta |K|}}(u)$.
Now, as in the proof of Theorem 6, we change the coordinate system
so that the Santal${\mbox{\'o}}$-point of K is at the origin and the positive
u-direction is the positive s-direction.
By construction of the floating body
$$\delta |K| \leq \int_{h_{K_{\delta |K|}}(u)}^{h_{K}(u)} g_{K,u}(s) ds$$
$$\leq v_{n-1} \int_{h_{K_{\delta |K|}}(u)}^{h_{K}(u)}
\biggl((\rho + \varepsilon)^2- (s+ \rho + \varepsilon
-h_{K}(u))^2\biggr)^{\frac{n-1}{2}}ds$$
$$\leq \frac{2^\frac{n+1}{2} v_{n-1}}{n+1} (\rho +\varepsilon)^\frac{n-1}{2}
\biggl(h_{K}(u)-h_{K_{\delta |K|}}(u)\biggr)^\frac{n+1}{2}.$$
Thus
$$\frac{1}{h_{K_{\delta |K|}}(u)^n} \geq \frac{1}{h_{K}(u)^n}\biggl(
1+ \frac {n (\delta|K|)^\frac{2}{n+1}}{c_n h_{K}(u)(\rho +
\varepsilon)^\frac{n-1}{n+1}}
\biggr)$$
and hence
\begin{equation}
\frac{c_n}{n (\delta|K|)^\frac{2}{n+1}} \biggl[ \frac{1}{h_{K_{\delta
|K|}}(u)^n}
-\frac{1}{h_{K}(u)^n}\biggr] \geq \frac{1}{h_{K}(u)^{n+1}(\rho
+\varepsilon)^\frac{n-1}{n+1}}.
\end{equation}
\vskip 3mm
On the other hand by the fact that for each $z \in
\partial K_{\delta |K|}$ there exists a tangent hyperplane that cuts off
exactly $\delta |K|$ of $K$ and by [S-W], Lemma 11, for $\varepsilon >0$
there exists $\delta_0$ such that for all $\delta \leq \delta _0$
$$\delta |K| \geq
(1-\varepsilon) v_{n-1} \int_{h_{K_{\delta |K|}}(u)+\varepsilon}^{h_{K}(u)}
\biggl((\rho - \varepsilon)^2- (s+ \rho - \varepsilon
-h_{K}(u))^2\biggr)^{\frac{n-1}{2}}ds$$
$$\geq \frac{(1-\varepsilon)2^\frac{n+1}{2} v_{n-1}}{n+1} (\rho
-\varepsilon)^\frac{n-1}{2}
\biggl(h_{K}(u)-h_{K_{\delta |K|}}(u)-\varepsilon\biggr)^\frac{n+1}{2}$$
$$\biggl(1-\frac{1}{\rho-\varepsilon}(h_{K}(u)-h_{K_{\delta
|K|}}(u)-\varepsilon)\biggr)^\frac{n-1}{2}$$
$$\geq \frac{(1-c \varepsilon)2^\frac{n+1}{2} v_{n-1}}{n+1} (\rho
-\varepsilon)^\frac{n-1}{2}
\biggl(h_{K}(u)-h_{K_{\delta |K|}}(u)\biggr)^\frac{n+1}{2},$$
for some constant c.
Thus for $\delta$ small enough
$$\frac{1}{h_{K_{\delta |K|}}(u)^n} \leq \frac{1}{h_{K}(u)^n}
\biggl\{
1+ \frac {n (\delta|K|)^\frac{2}{n+1}(1+ dA(\delta))}{c_n (1-c
\varepsilon)^\frac{2}{n+1}
h_{K}(u)(\rho -
\varepsilon)^\frac{n-1}{n+1}}$$
$$\biggl(1 +\frac{1}{n}
\left (
\begin{array}{c}
 n \\ 2
\end{array}
\right)
 A(\delta)(1+d A(\delta)) +
\frac{1}{n}
\left (
\begin{array}{c}
 n \\ 3
\end{array}
\right)
( A(\delta))^2(1+d A(\delta))^2  +\cdots
$$
$$\cdots +
\frac{1}{n} ( A(\delta)^n(1+d A(\delta))^n
\biggr)\biggr\},$$

where d is a constant and $A(\delta )=  \frac {n
(\delta|K|)^\frac{2}{n+1}}{c_n (1-c
\varepsilon)^\frac{2}{n+1} h_{K}(u)(\rho -
\varepsilon)^\frac{n-1}{n+1}}$.
\newline
Thus
$$ \frac{c_n}{n (\delta|K|)^\frac{2}{n+1}} \biggl[ \frac{1}{h_{K_{\delta
|K|}}(u)^n}
-\frac{1}{h_{K}(u)^n}\biggr] \leq
\frac{(1+ d A(\delta))}{(1-c \varepsilon)^\frac{2}{n-1}h_{K}(u)^{n+1}(\rho
-\varepsilon)^\frac{n-1}{n+1}}$$
$$\biggl\{
1 +\frac{1}{n}
\left (
\begin{array}{c}
 n \\ 2
\end{array}
\right)
 A(\delta)(1+d A(\delta)) +
\frac{1}{n}
\left (
\begin{array}{c}
 n \\ 3
\end{array}
\right)
( A(\delta))^2(1+d A(\delta))^2  +\cdots
$$
\begin{equation}
\cdots +
\frac{1}{n} ( A(\delta)^n(1+d A(\delta))^n
\biggr\}.
\end{equation}
\vskip 5mm
(12) and (13) show that for $u \in S^{n-1}$

$$\mbox{lim}_{\delta \rightarrow 0}
\frac{c_n}{n (\delta|K|)^\frac{2}{n+1}} \biggl[
\frac{1}{h_{K_{\delta |K|}}(u)^n} -\frac{1}{h_{K}(u)^n}\biggr] = \frac{\rho
(u)^{-\frac{n-1}{n+1}}}{h_{K}(u)^{n+1}}.$$
Therefore
$$\mbox{lim}_{\delta \rightarrow 0}
c_{n}\frac{|(K_{\delta |K|})^0| -|K^0|}{(\delta |K|)^{\frac{2}{n+1}}}=$$
$$\mbox{lim}_{\delta \rightarrow 0}
\int_{S^{n-1}}\frac{c_n}{n (\delta|K|)^\frac{2}{n+1}} \biggl[
\frac{1}{h_{K_{\delta |K|}}(u)^n} -\frac{1}{h_{K}(u)^n}\biggr]
d\sigma(u)=$$
$$\int_{S^{n-1}}\mbox{lim}_{\delta \rightarrow 0}
\frac{c_n}{n (\delta|K|)^\frac{2}{n+1}} \biggl[
\frac{1}{h_{K_{\delta |K|}}(u)^n} -\frac{1}{h_{K}(u)^n}\biggr]
d\sigma(u)$$
$$=
\int_{S^{n-1}} \frac{d\sigma(u)}{
f_{K}(u)^{\frac{1}{n+1}}h_{K}(u)^{n+1}}.$$
\newpage
We still have to justify that we can interchange integration and limit.
This follows  from
Lebesgue's Theorem and the
\newline
\underline{Claim}
\newline
for every $u \in S^{n-1}$
$$\frac{1}{n (\delta|K|)^\frac{2}{n+1}} \biggl[
\frac{1}{h_{K_{\delta |K|}}(u)^n} -\frac{1}{h_{K}(u)^n}\biggr]  \leq g(u),$$
where g is a function independent of $\delta$ and integrable on $S^{n-1}$.
\newline
\underline{Proof of the Claim}
\par
For $u \in S^{n-1}$ let $ x \in \partial K$ be such that $N(x)=u$.
Moreover we can suppose that $0 \in \mbox{int} (K)$ and choose  $\alpha >
0$ such
that
\begin{equation}
B(0,\alpha) \subseteq K \subseteq
B(0,\frac{1}{\alpha}).
\end{equation}
Let again $R_0 = \mbox{min}_{x \in \partial K,1 \leq i \leq n-1}R_i(x)$,
where $R_i(x)$
is the i-th principal radius of curvature at $x \in \partial K$.
$R_0 > 0$  (see [L]).
\newline
By the Blaschke Rolling Theorem (see [L]), we have for all $x \in \partial K$
\begin{equation}
B(x-R_0N(x), R_0) \subseteq K.
\end{equation}
By  [L] there exists $\delta_0$ such that for all $\delta < \delta_0$
$\partial K_{\delta |K|}$ is $C^2$.
\newline
We put
$\delta_1= \mbox{min}\{\delta_0, (\frac{3}{4}
\alpha)^{\frac{n+1}{2}}\frac{c_nR_0^{\frac{n-1}{2}}}{|K|}\}$.
Then we have for all $\delta < \delta_1$ that $\partial K_{\delta |K|}$ is
$C^2$.
Consequently for $u \in S^{n-1}$ there is  $z \in \partial K_{\delta |K|}$
such that $N(z)=u$ and the tangent-hyperplane  to $K_{\delta |K|}$ in $z$
orthogonal to $u$ cuts off exactly $\delta   |K|$ from $K$.
\newline
Now we distinguish two cases.
\newline
a)  $|| x-z || <\frac{x}{ ||x ||},u>  \leq R_0$.
\newline
Because of (15), $ \delta  |K| $ can be estimated from below by the volume
of the cap of
height $|| x-z || <\frac{x}{ ||x ||},u>= h_{K}(u)-h_{K_{\delta |K|}}(u)$ of
a Euclidean ball with
radius
$R_0$. Thus
$$\delta  |K| \geq \frac{2v_{n-1}}{n+1}R_0^{\frac{n-1}{2}}
\biggl(h_{K}(u)-h_{K_{\delta |K|}}(u) \biggr)^{\frac{n+1}{2}},$$
and hence with (14) and the choice of $\delta_1$
$$\frac{1}{(h_{K_{\delta |K|}}(u))^n} \leq \frac{1+c n (\delta |K|
)^{\frac{2}{n+1}}}{(h_{K}(u))^n},$$
where $c$ is a constant independent of $\delta $ and $u$.
Therefore we get in this case with a constant $d$
$$\frac{1}{n (\delta|K|)^\frac{2}{n+1}} \biggl[
\frac{1}{h_{K_{\delta |K|}}(u)^n} -\frac{1}{h_{K}(u)^n}\biggr]  \leq \frac{d
}{R_0^\frac{n-1}{n+1} h_K(u)^{n+1}}$$
and the latter is an integrable function on $S^{n-1}$.
\vskip 3mm
b)  $|| x-z || <\frac{x}{ ||x ||},u>  > R_0$.
\newline
Let $H_x$ be the hyperplane through $0$ orthogonal to the vector $x$
and let $C$ be the cone $C=\mbox{co}[x,H_x \cap B(0,\alpha)]$.
Then the tangent-hyperplane to $\partial K_{\delta |K|}$ through $z$
cuts off more from $K$ than it does from $C$ and therefore
$$\delta  |K| \geq \frac{v_{n-1}\alpha ^{n-1}|| x-z ||^n}{n || x|| ^{n-1}},$$
as the volume that a hyperplane cuts off the cone $C$ is minimal if the
hyperplane is parallel to the base of the cone.
Thus
$$\delta  |K| \geq \frac{1}{n}v_{n-1} \alpha ^{2(n-1)}
\biggl(h_{K}(u)-h_{K_{\delta |K|}}(u) \biggr)^{n}$$
and hence with a constant $k$ (independent of $u$ and $\delta$)
$$\frac{1}{n (\delta|K|)^\frac{2}{n+1}} \biggl[
\frac{1}{h_{K_{\delta |K|}}(u)^n} -\frac{1}{h_{K}(u)^n}\biggr]  \leq \frac{k
}{ (\delta|K|)^{\frac{n-1}{n(n+1)}}(h_K(u))^{n}} .$$
By assumption b) $\delta  |K| \geq \frac{1}{2}R_0^nv_{n}$ and thus the
latter can be
estimated by
$$\leq k (\frac{2}{v_n})^{\frac{n-1}{n(n+1)}}\frac{1}{R_0^\frac{n-1}{n+1}}
\frac{1}{h_K(u)^{n}},$$
which is integrable on $S^{n-1}$.
\vskip 3mm
This finishes the proof of the Claim and thus of the Theorem.

\newpage
\medskip\noindent Mathieu Meyer
\newline
Universit\'{e} de Marne-la-Vall\'{e}e
\newline
Equipe d'Analyse
\newline
93166 Noisy-le-Grand Cedex, France
\newline
e-mail:meyer@math.univ.mlv.fr
\vskip 3mm
\medskip\noindent Elisabeth Werner
\newline
Department of Mathematics
\newline
Case Western Reserve University
\newline
Cleveland, Ohio 44106, U.S.A.
\newline
e-mail: emw2@po.cwru.edu
\newline
and
\newline
Universit\'{e} de Lille 1
\newline
Ufr de Mathematique
\newline
59655 Villeneuve d'Ascq, France

\end{document}